\documentclass{amsart} 
\usepackage{amsfonts,amssymb,amsmath,euscript}
    \newtheorem{thm}{Theorem}                     [section]
    \newtheorem{thm*}{Theorem}
    \newtheorem{prop}[thm]{Proposition}
    \newtheorem{lemma}[thm]{Lemma}
    \newtheorem{cor}[thm]{Corollary}

    \newtheorem{lemma*}{Lemma}    

    \newtheorem{assump}[thm]{Assumption}

    \newtheorem{defn}[thm]{Definition}                 
    
    \newtheorem{rems*}{Remark}   

\renewenvironment{proof}{{\it Proof.}}{$\Box$\par} 

\newcommand{\ndef}{\newcommand*}
\def\rndef{\renewcommand}

\ndef{\myaddress}[1]{\begin{center} \it\small #1 \end{center}}

\ndef{\clA}{{\mathcal A}} \ndef{\rmA}{{\mathrm A}} \ndef{\mbA}{{\mathbb A}} \ndef{\bfA}{{\mathbf A}} \ndef{\euA}{{\EuScript A}} \ndef{\frA}{{\mathfrak A}}
\ndef{\clB}{{\mathcal B}} \ndef{\rmB}{{\mathrm B}} \ndef{\mbB}{{\mathbb B}} \ndef{\bfB}{{\mathbf B}} \ndef{\euB}{{\EuScript B}} \ndef{\frB}{{\mathfrak B}}
\ndef{\clC}{{\mathcal C}} \ndef{\rmC}{{\mathrm C}} \ndef{\mbC}{{\mathbb C}} \ndef{\bfC}{{\mathbf C}} \ndef{\euC}{{\EuScript C}} \ndef{\frC}{{\mathfrak C}}
\ndef{\clD}{{\mathcal D}} \ndef{\rmD}{{\mathrm D}} \ndef{\mbD}{{\mathbb D}} \ndef{\bfD}{{\mathbf D}} \ndef{\euD}{{\EuScript D}} \ndef{\frD}{{\mathfrak D}}
\ndef{\clE}{{\mathcal E}} \ndef{\rmE}{{\mathrm E}} \ndef{\mbE}{{\mathbb E}} \ndef{\bfE}{{\mathbf E}} \ndef{\euE}{{\EuScript E}} \ndef{\frE}{{\mathfrak E}}
\ndef{\clF}{{\mathcal F}} \ndef{\rmF}{{\mathrm F}} \ndef{\mbF}{{\mathbb F}} \ndef{\bfF}{{\mathbf F}} \ndef{\euF}{{\EuScript F}} \ndef{\frF}{{\mathfrak F}}
\ndef{\clG}{{\mathcal G}} \ndef{\rmG}{{\mathrm G}} \ndef{\mbG}{{\mathbb G}} \ndef{\bfG}{{\mathbf G}} \ndef{\euG}{{\EuScript G}} \ndef{\frG}{{\mathfrak G}}
\ndef{\clH}{{\mathcal H}} \ndef{\rmH}{{\mathrm H}} \ndef{\mbH}{{\mathbb H}} \ndef{\bfH}{{\mathbf H}} \ndef{\euH}{{\EuScript H}} \ndef{\frH}{{\mathfrak H}}
\ndef{\clI}{{\mathcal I}} \ndef{\rmI}{{\mathrm I}} \ndef{\mbI}{{\mathbb I}} \ndef{\bfI}{{\mathbf I}} \ndef{\euI}{{\EuScript I}} \ndef{\frI}{{\mathfrak I}}
\ndef{\clJ}{{\mathcal J}} \ndef{\rmJ}{{\mathrm J}} \ndef{\mbJ}{{\mathbb J}} \ndef{\bfJ}{{\mathbf J}} \ndef{\euJ}{{\EuScript J}} \ndef{\frJ}{{\mathfrak J}}
\ndef{\clK}{{\mathcal K}} \ndef{\rmK}{{\mathrm K}} \ndef{\mbK}{{\mathbb K}} \ndef{\bfK}{{\mathbf K}} \ndef{\euK}{{\EuScript K}} \ndef{\frK}{{\mathfrak K}}
\ndef{\clL}{{\mathcal L}} \ndef{\rmL}{{\mathrm L}} \ndef{\mbL}{{\mathbb L}} \ndef{\bfL}{{\mathbf L}} \ndef{\euL}{{\EuScript L}} \ndef{\frL}{{\mathfrak L}}
\ndef{\clM}{{\mathcal M}} \ndef{\rmM}{{\mathrm M}} \ndef{\mbM}{{\mathbb M}} \ndef{\bfM}{{\mathbf M}} \ndef{\euM}{{\EuScript M}} \ndef{\frM}{{\mathfrak M}}
\ndef{\clN}{{\mathcal N}} \ndef{\rmN}{{\mathrm N}} \ndef{\mbN}{{\mathbb N}} \ndef{\bfN}{{\mathbf N}} \ndef{\euN}{{\EuScript N}} \ndef{\frN}{{\mathfrak N}}
\ndef{\clO}{{\mathcal O}} \ndef{\rmO}{{\mathrm O}} \ndef{\mbO}{{\mathbb O}} \ndef{\bfO}{{\mathbf O}} \ndef{\euO}{{\EuScript O}} \ndef{\frO}{{\mathfrak O}}
\ndef{\clP}{{\mathcal P}} \ndef{\rmP}{{\mathrm P}} \ndef{\mbP}{{\mathbb P}} \ndef{\bfP}{{\mathbf P}} \ndef{\euP}{{\EuScript P}} \ndef{\frP}{{\mathfrak P}}
\ndef{\clQ}{{\mathcal Q}} \ndef{\rmQ}{{\mathrm Q}} \ndef{\mbQ}{{\mathbb Q}} \ndef{\bfQ}{{\mathbf Q}} \ndef{\euQ}{{\EuScript Q}} \ndef{\frQ}{{\mathfrak Q}}
\ndef{\clR}{{\mathcal R}} \ndef{\rmR}{{\mathrm R}} \ndef{\mbR}{{\mathbb R}} \ndef{\bfR}{{\mathbf R}} \ndef{\euR}{{\EuScript R}} \ndef{\frR}{{\mathfrak R}}
\ndef{\clS}{{\mathcal S}} \ndef{\rmS}{{\mathrm S}} \ndef{\mbS}{{\mathbb S}} \ndef{\bfS}{{\mathbf S}} \ndef{\euS}{{\EuScript S}} \ndef{\frS}{{\mathfrak S}}
\ndef{\clT}{{\mathcal T}} \ndef{\rmT}{{\mathrm T}} \ndef{\mbT}{{\mathbb T}} \ndef{\bfT}{{\mathbf T}} \ndef{\euT}{{\EuScript T}} \ndef{\frT}{{\mathfrak T}}
\ndef{\clU}{{\mathcal U}} \ndef{\rmU}{{\mathrm U}} \ndef{\mbU}{{\mathbb U}} \ndef{\bfU}{{\mathbf U}} \ndef{\euU}{{\EuScript U}} \ndef{\frU}{{\mathfrak U}}
\ndef{\clV}{{\mathcal V}} \ndef{\rmV}{{\mathrm V}} \ndef{\mbV}{{\mathbb V}} \ndef{\bfV}{{\mathbf V}} \ndef{\euV}{{\EuScript V}} \ndef{\frV}{{\mathfrak V}}
\ndef{\clW}{{\mathcal W}} \ndef{\rmW}{{\mathrm W}} \ndef{\mbW}{{\mathbb W}} \ndef{\bfW}{{\mathbf W}} \ndef{\euW}{{\EuScript W}} \ndef{\frW}{{\mathfrak W}}
\ndef{\clX}{{\mathcal X}} \ndef{\rmX}{{\mathrm X}} \ndef{\mbX}{{\mathbb X}} \ndef{\bfX}{{\mathbf X}} \ndef{\euX}{{\EuScript X}} \ndef{\frX}{{\mathfrak X}}
\ndef{\clY}{{\mathcal Y}} \ndef{\rmY}{{\mathrm Y}} \ndef{\mbY}{{\mathbb Y}} \ndef{\bfY}{{\mathbf Y}} \ndef{\euY}{{\EuScript Y}} \ndef{\frY}{{\mathfrak Y}}
\ndef{\clZ}{{\mathcal Z}} \ndef{\rmZ}{{\mathrm Z}} \ndef{\mbZ}{{\mathbb Z}} \ndef{\bfZ}{{\mathbf Z}} \ndef{\euZ}{{\EuScript Z}} \ndef{\frZ}{{\mathfrak Z}}

\ndef{\tA}{{\widetilde A}} \ndef{\tcA}{{\widetilde\clA}} \ndef{\ttcA}{\widetilde{\tcA}} \ndef{\sfA}{{\textsf A}} \ndef{\ttA}{\widetilde{\tA}} \ndef{\dzA}{{A^\sharp}}
\ndef{\tB}{{\widetilde B}} \ndef{\tcB}{{\widetilde\clB}} \ndef{\ttcB}{\widetilde{\tcB}} \ndef{\sfB}{{\textsf B}} \ndef{\ttB}{\widetilde{\tB}} \ndef{\dzB}{{B^\sharp}}
\ndef{\tC}{{\widetilde C}} \ndef{\tcC}{{\widetilde\clC}} \ndef{\ttcC}{\widetilde{\tcC}} \ndef{\sfC}{{\textsf C}} \ndef{\ttC}{\widetilde{\tC}} \ndef{\dzC}{{C^\sharp}}
\ndef{\tD}{{\widetilde D}} \ndef{\tcD}{{\widetilde\clD}} \ndef{\ttcD}{\widetilde{\tcD}} \ndef{\sfD}{{\textsf D}} \ndef{\ttD}{\widetilde{\tD}} \ndef{\dzD}{{D^\sharp}}
\ndef{\tE}{{\widetilde E}} \ndef{\tcE}{{\widetilde\clE}} \ndef{\ttcE}{\widetilde{\tcE}} \ndef{\sfE}{{\textsf E}} \ndef{\ttE}{\widetilde{\tE}} \ndef{\dzE}{{E^\sharp}}
\ndef{\tF}{{\widetilde F}} \ndef{\tcF}{{\widetilde\clF}} \ndef{\ttcF}{\widetilde{\tcF}} \ndef{\sfF}{{\textsf F}} \ndef{\ttF}{\widetilde{\tF}} \ndef{\dzF}{{F^\sharp}}
\ndef{\tG}{{\widetilde G}} \ndef{\tcG}{{\widetilde\clG}} \ndef{\ttcG}{\widetilde{\tcG}} \ndef{\sfG}{{\textsf G}} \ndef{\ttG}{\widetilde{\tG}} \ndef{\dzG}{{G^\sharp}}
\ndef{\tH}{{\widetilde H}} \ndef{\tcH}{{\widetilde\clH}} \ndef{\ttcH}{\widetilde{\tcH}} \ndef{\sfH}{{\textsf H}} \ndef{\ttH}{\widetilde{\tH}} \ndef{\dzH}{{H^\sharp}}
\ndef{\tI}{{\widetilde I}} \ndef{\tcI}{{\widetilde\clI}} \ndef{\ttcI}{\widetilde{\tcI}} \ndef{\sfI}{{\textsf I}} \ndef{\ttI}{\widetilde{\tI}} \ndef{\dzI}{{I^\sharp}}
\ndef{\tJ}{{\widetilde J}} \ndef{\tcJ}{{\widetilde\clJ}} \ndef{\ttcJ}{\widetilde{\tcJ}} \ndef{\sfJ}{{\textsf J}} \ndef{\ttJ}{\widetilde{\tJ}} \ndef{\dzJ}{{J^\sharp}}
\ndef{\tK}{{\widetilde K}} \ndef{\tcK}{{\widetilde\clK}} \ndef{\ttcK}{\widetilde{\tcK}} \ndef{\sfK}{{\textsf K}} \ndef{\ttK}{\widetilde{\tK}} \ndef{\dzK}{{K^\sharp}}
\ndef{\tL}{{\widetilde L}} \ndef{\tcL}{{\widetilde\clL}} \ndef{\ttcL}{\widetilde{\tcL}} \ndef{\sfL}{{\textsf L}} \ndef{\ttL}{\widetilde{\tL}} \ndef{\dzL}{{L^\sharp}}
\ndef{\tM}{{\widetilde M}} \ndef{\tcM}{{\widetilde\clM}} \ndef{\ttcM}{\widetilde{\tcM}} \ndef{\sfM}{{\textsf M}} \ndef{\ttM}{\widetilde{\tM}} \ndef{\dzM}{{M^\sharp}}
\ndef{\tN}{{\widetilde N}} \ndef{\tcN}{{\widetilde\clN}} \ndef{\ttcN}{\widetilde{\tcN}} \ndef{\sfN}{{\textsf N}} \ndef{\ttN}{\widetilde{\tN}} \ndef{\dzN}{{N^\sharp}}
\ndef{\tO}{{\widetilde O}} \ndef{\tcO}{{\widetilde\clO}} \ndef{\ttcO}{\widetilde{\tcO}} \ndef{\sfO}{{\textsf O}} \ndef{\ttO}{\widetilde{\tO}} \ndef{\dzO}{{O^\sharp}}
\ndef{\tP}{{\widetilde P}} \ndef{\tcP}{{\widetilde\clP}} \ndef{\ttcP}{\widetilde{\tcP}} \ndef{\sfP}{{\textsf P}} \ndef{\ttP}{\widetilde{\tP}} \ndef{\dzP}{{P^\sharp}}
\ndef{\tQ}{{\widetilde Q}} \ndef{\tcQ}{{\widetilde\clQ}} \ndef{\ttcQ}{\widetilde{\tcQ}} \ndef{\sfQ}{{\textsf Q}} \ndef{\ttQ}{\widetilde{\tQ}} \ndef{\dzQ}{{Q^\sharp}}
\ndef{\tR}{{\widetilde R}} \ndef{\tcR}{{\widetilde\clR}} \ndef{\ttcR}{\widetilde{\tcR}} \ndef{\sfR}{{\textsf R}} \ndef{\ttR}{\widetilde{\tR}} \ndef{\dzR}{{R^\sharp}}
\ndef{\tS}{{\widetilde S}} \ndef{\tcS}{{\widetilde\clS}} \ndef{\ttcS}{\widetilde{\tcS}} \ndef{\sfS}{{\textsf S}} \ndef{\ttS}{\widetilde{\tS}} \ndef{\dzS}{{S^\sharp}}
\ndef{\tT}{{\widetilde T}} \ndef{\tcT}{{\widetilde\clT}} \ndef{\ttcT}{\widetilde{\tcT}} \ndef{\sfT}{{\textsf T}} \ndef{\ttT}{\widetilde{\tT}} \ndef{\dzT}{{T^\sharp}}
\ndef{\tU}{{\widetilde U}} \ndef{\tcU}{{\widetilde\clU}} \ndef{\ttcU}{\widetilde{\tcU}} \ndef{\sfU}{{\textsf U}} \ndef{\ttU}{\widetilde{\tU}} \ndef{\dzU}{{U^\sharp}}
\ndef{\tV}{{\widetilde V}} \ndef{\tcV}{{\widetilde\clV}} \ndef{\ttcV}{\widetilde{\tcV}} \ndef{\sfV}{{\textsf V}} \ndef{\ttV}{\widetilde{\tV}} \ndef{\dzV}{{V^\sharp}}
\ndef{\tW}{{\widetilde W}} \ndef{\tcW}{{\widetilde\clW}} \ndef{\ttcW}{\widetilde{\tcW}} \ndef{\sfW}{{\textsf W}} \ndef{\ttW}{\widetilde{\tW}} \ndef{\dzW}{{W^\sharp}}
\ndef{\tX}{{\widetilde X}} \ndef{\tcX}{{\widetilde\clX}} \ndef{\ttcX}{\widetilde{\tcX}} \ndef{\sfX}{{\textsf X}} \ndef{\ttX}{\widetilde{\tX}} \ndef{\dzX}{{X^\sharp}}
\ndef{\tY}{{\widetilde Y}} \ndef{\tcY}{{\widetilde\clY}} \ndef{\ttcY}{\widetilde{\tcY}} \ndef{\sfY}{{\textsf Y}} \ndef{\ttY}{\widetilde{\tY}} \ndef{\dzY}{{Y^\sharp}}
\ndef{\tZ}{{\widetilde Z}} \ndef{\tcZ}{{\widetilde\clZ}} \ndef{\ttcZ}{\widetilde{\tcZ}} \ndef{\sfZ}{{\textsf Z}} \ndef{\ttZ}{\widetilde{\tZ}} \ndef{\dzZ}{{Z^\sharp}}

\ndef{\bfc}{{\bf c}}

  \ndef{\eps}{\varepsilon}

\let\leq\leqslant

\ndef{\lims}[1]{\lim\limits_{#1}}
\ndef{\sums}[1]{\sum\limits_{#1}}
\ndef{\ints}[1]{\int\limits_{#1}}
\ndef{\sups}[1]{\sup\limits_{#1}}
\ndef{\liminfty}[1]{\lims{#1\to\infty}}
\ndef{\suminf}[1]{\sums{#1=1}^\infty}

\ndef{\limo}[1]{\omega\mbox{-}\!\!\!\lims{#1\to\infty}}          
\ndef{\limL}[1]{\rmL\mbox{-}\!\!\!\lims{#1\to\infty}}            
\ndef{\limLOne}[1]{\clL_1\mbox{-}\!\lims{#1}}
\ndef{\tildelimo}[1]{\tilde\omega\mbox{-}\!\!\!\lims{#1\to\infty}}
\ndef{\slim}{\mathrm{s}\mbox{-}\!\!\lim}          
\ndef{\wlim}{\mathrm{w}\mbox{-}\!\lim}          

\ndef{\Aut}{\operatorname{Aut}}      
\ndef{\Ch}{\operatorname{ch}}        
\ndef{\End}{\operatorname{End}}      
\ndef{\Hom}{\operatorname{Hom}}      
\ndef{\Ker}{\operatorname{Ker}}      
\ndef{\Log}{\operatorname{Log}}      
\ndef{\OP}{\operatorname{OP}}        
\ndef{\Op}{\operatorname{Op}}        
\ndef{\Symb}{\operatorname{Symb}}    
\ndef{\Tr}{\operatorname{Tr}}        
\ndef{\Wres}{\operatorname{Wres}}    
\ndef{\cl}{\operatorname{cl}}        
\ndef{\com}{\operatorname{com}}
\ndef{\const}{\operatorname{const}}  
\ndef{\conv}{\operatorname{conv}}    
\rndef{\det}{\operatorname{det}}     

\ndef{\detFK}[1]{\Delta\brs{#1}} 
\ndef{\detFKrel}[2]{\Delta_{#2}\brs{#1}} 

\ndef{\adj}{\operatorname{adj}}    
\ndef{\diag}{\operatorname{diag}}    
\ndef{\dist}{\operatorname{dist}}    
\ndef{\dom}{\operatorname{dom}}      
\ndef{\ec}{\operatorname{ec}}        
\ndef{\id}{1}                        
\ndef{\ind}{\operatorname{ind}}      
\ndef{\mydeg}{\operatorname{deg}}    
\ndef{\op}{\operatorname{op}}
\ndef{\rank}{\operatorname{rank}}
\ndef{\res}{\operatorname{res}}      
\ndef{\rng}{\operatorname{ran}}      
\ndef{\sflow}{\operatorname{sf}}     
\ndef{\isf}{\operatorname{isf}}      
\ndef{\sign}{\operatorname{sign}}    
\ndef{\sgn}{\operatorname{sgn}}      
\ndef{\sing}{\operatorname{sing}}    
\ndef{\supp}{\operatorname{supp}}    
\ndef{\tr}{\operatorname{tr}}        
\ndef{\var}{\operatorname{var}}      
\ndef{\vol}{\operatorname{vol}}      
\ndef{\wn}{\operatorname{wn}}        
\ndef{\wres}{\operatorname{wres}}    
\rndef{\Im}{\operatorname{Im}}       
\rndef{\Re}{\operatorname{Re}}       

\ndef{\prng}[1]{\mathrm R_{#1}} 
\ndef{\pker}[1]{\mathrm N_{#1}} 
\ndef{\rprng}[2]{\mathrm R_{#1}^{#2}}           
\ndef{\rpker}[2]{\mathrm N_{#1}^{#2}}           
\ndef{\rsupp}[1]{\supp_r(#1)}
\ndef{\lsupp}[1]{\supp_l(#1)}
\ndef{\rslv}[1]{R_z(#1)}      
\ndef{\HH}{H}                 
\ndef{\tHH}{\tilde \HH}       
\ndef{\VV}{V}                 
\ndef{\Rz}{R_z}               
\ndef{\tRz}{\tR_z}            
\ndef{\psif}[1]{#1^{[1]}} 
\ndef{\CPlus}[1]{W_{#1}(\mbR)}
\ndef{\bndl}{\xi}                         
\ndef{\bndlA}{\eta}                       
\ndef{\GlueMap}{\varphi}                  
\ndef{\ChartMap}{h}                       

\ndef{\hilb}{\clH}                     
\ndef{\hilba}{\clH^{(a)}}                    
   \ndef{\hilbasargument}{(\hilb)} 
\ndef{\LpH}[1]{\clL^{#1}\hilbasargument}       
\ndef{\saLpH}[1]{\clL_{sa}^{#1}\hilbasargument}       
\ndef{\clBH}{\clB\hilbasargument}              
\ndef{\ubBH}{\clB_1\hilbasargument}            
\ndef{\clCH}{\clC\hilbasargument}              
\ndef{\clKH}{\clK\hilbasargument}              
\ndef{\clFH}{\clF\hilbasargument}              
\ndef{\clUH}{\clU\hilbasargument}              
\ndef{\clCFH}{{\clC\clF}\hilbasargument}       
\ndef{\saBH}{\clB_{sa}\hilbasargument}         
\ndef{\saCH}{\clC_{sa}\hilbasargument}         
\ndef{\saFH}{\clF_{sa}\hilbasargument}         
\ndef{\saKH}{\clK_{sa}\hilbasargument}         
\ndef{\saCFH}{\clC\clF_{sa}\hilbasargument}    
\ndef{\clUFH}{\clU\clF\hilbasargument}         
\ndef{\Uinj}{\clU_{inj}\hilbasargument}        
\ndef{\UFinj}{\clU\clF_{inj}\hilbasargument}   

\ndef{\spproj}[2]{E^{#1}_{#2}}                      
\ndef{\spprojb}[2]{E^{#2}_{#1}}                     

\ndef{\LpN}[1]{\clL^{#1}(\clN,\tau)}     
\ndef{\saLpN}[1]{\clL^{#1}_{sa}(\clN,\tau)} 
\ndef{\rLpN}[1]{L^{#1}(\clN,\tau)}       
\ndef{\clAND}{(\clA,\clN,D)}             
\ndef{\clBA}{{\clB(\clA)}}
\ndef{\saKN}{{\clK_{sa}(\clN,\tau)}}          
\ndef{\clKN}{{\clK(\clN,\tau)}}          
\ndef{\clKtN}{{\clK(\tilde\clN,\tau)}}   
\ndef{\clFN}{{\clF(\clN,\tau)}}          
\ndef{\saFN}{{\clF_{sa}(\clN,\tau)}}     
\ndef{\clPN}{\clP(\clN)}                 
\ndef{\clQN}{\clQ(\clN,\tau)}            
\ndef{\infPN}{{\clP_\tau^\infty(\clN)}}  
\ndef{\clOF}[2]{\clF_{#1\mbox{-}#2}(\clN,\tau)}         
\ndef{\oind}[2]{{\rm \tau\mbox{-}ind}_{#1\mbox{-}#2}}   
\ndef{\tind}{\tau\mbox{-}\ind}                  
\ndef{\DInd}{\ind_{\clD,\tau}}           
\ndef{\BF}{Breuer-Fredholm}              
\ndef{\skewfred}[2]{$(#1\cdot #2)$ $\tau$\tire Fredholm}   
\ndef{\affl}{\eta}                       
\ndef{\vNa}{von Neumann algebra}         
\ndef{\nsf}{faithful normal semifinite } 
\ndef{\taubrs}[1]{\tau\brackets{#1}}     
\ndef{\sqbrs}[1]{[#1]}        
\ndef{\Sqbrs}[1]{\big[#1\big]}        
\ndef{\SqBrs}[1]{\Big[#1\Big]}        

\ndef{\domd}{\bigcap\limits_{n\ge 0} \dom\;\delta^n}         
\ndef{\DiffOP}{{\rm \clD}}
\ndef{\ADA}{\clA \cup [\clD,\clA]}
\ndef{\DixIdeal}[1]{\LpH{#1,\infty}}               
\ndef{\dixideal}{\ell^{1,\infty}}                  
\ndef{\WDixIdeal}{\LpH{1,\mathrm w}}               
\ndef{\DixIdealPos}[1]{\DixIdeal{#1}_+}            
\ndef{\DixIdealN}[1]{\LpN{#1,\infty}}              
\ndef{\DixIdealNPar}[2]{\clL^{#1,\infty}_{#2}(\clN,\tau)}    
\ndef{\DixIdealNPos}[1]{\LpN{#1,\infty}_+}                   
\ndef{\TrD}{\Tr_\omega}                                      
\ndef{\tauD}{{\tau_\omega}}                                  
\ndef{\ILog}{\frac 1{\log(1+t)}}
\ndef{\ILogN}{\frac 1{\log(1+N)}}
\ndef{\DixNorm}[1]{\norm{#1}_{(1,\infty)}}                   
\ndef{\DixInt}[1]{\ints 0^t \mu_s(#1)\,ds}
\ndef{\DixIntL}[1]{\ints 0^{\lambda_{1/t}(#1)}\mu_s(#1)\,ds}
    \ndef{\SmallIdeal}{{\clL^{1, \mathrm w}}}
    \ndef{\SmallIdealMeas}{{\clL^{1, \mathrm w}_m}}
    \ndef{\DixIntII}[1]{\ints 0^t \mu_s(#1)\,ds}
    \ndef{\DixIntf}[1]{f_t(#1)}
    \ndef{\DixIntg}[1]{g_t(#1)}

\ndef{\lpi}{\clL^{1,\pi}(\clN,\tau)}
\ndef{\IIinfty}{$\mathrm{II}_\infty$\ }

\ndef{\fourier}[1]{\clF(#1)}          
\ndef{\HaarMeasBohrs}{\nu}            
\ndef{\BrownsMeas}{\mu}               
\ndef{\BohrCont}[1]{\tilde{#1}}       
\ndef{\APMean}{{M}}                   
\ndef{\CDSS}{{\clA_B}}                
\ndef{\matr}{{\rm Mat}}               
\ndef{\seque}[1]{\ensuremath{\{#1_j\}_{j=1}^\infty}}    
\ndef{\sequen}[2]{\ensuremath{\{#1_#2\}_{#2=1}^\infty}}    
\ndef{\Seque}[1]{\ensuremath{\left(#1_0,#1_1,#1_2,\dots\right)}}    
\ndef{\Cesaro}{H}                           
\ndef{\CesaroRPlus}{M}                      
\ndef{\Dilation}{D}                         
\ndef{\Shift}{T}                            

\ndef{\norm}[1]{\left\Vert#1\right\Vert}    
\ndef{\TrNorm}[1]{\norm{#1}_1}              
\ndef{\HSNorm}[1]{\norm{#1}_2}              
\ndef{\InftyNorm}[1]{\norm{#1}_\infty}      
\ndef{\normQN}[1]{\norm{#1}_{\clQN}}        
\ndef{\clLnorm}[1]{\norm{#1}_{1,\infty}}    

\ndef{\ccurve}{\gamma}                      

\ndef{\abs}[1]{\left\lvert#1\right\rvert}   
\ndef{\set}[1]{\left\{#1\right\}}           
\ndef{\brackets}[1]{\left(#1\right)}        
\ndef{\brs}[1]{\brackets{#1}}               
\ndef{\Brs}[1]{\big(#1\big)}                
\ndef{\BRS}[1]{\Big(#1\Big)}                
\ndef{\scal}[2]{\left\la #1,#2\right\ra}               
\ndef{\precprec}{\prec\!\!\!\prec}
\ndef{\qeq}{\stackrel?=}
\ndef{\spectrum}[1]{\sigma_{#1}} 
\ndef{\spectruma}[1]{\sigma^{(a)}_{#1}} 
\ndef{\numrange}[1]{\mathrm{W}(#1)}                         
\rndef{\emptyset}{\varnothing}                              
\ndef{\csupp}{c}                           
\ndef{\closure}[1]{\overline{#1}}
\ndef{\linspan}[1]{\mathrm{span}\ {#1}}
\ndef{\bddborel}[1]{B(#1)}                 
\ndef{\charfunc}{\chi}
\rndef{\ln}{\log}
\ndef{\FrDer}{\euD}                        
\ndef{\LieDer}[1]{\pounds_{#1}\,}          
\ndef{\dds}{\left.\frac d{ds} \right|_{s = 0}}
\ndef{\ortcmp}[1]{#1^{\scriptscriptstyle \perp}}            
\ndef{\Laplace}{\Delta}                    

\ndef{\matrPQ}[3]
{
    \left(
      \begin{array}{cc}
        #1_{11} & #1_{12} \\
        #1_{21} & #1_{22}
      \end{array}
    \right)_{[#2,#3]}
}

\newcounter{margcomcount}
\ndef{\margcom}[1]{\marginpar{\bf \small #1} \addtocounter{margcomcount}{1}}

\newcounter{margproof}
\ndef{\margproof}{\marginpar{\bf \small Proof} \addtocounter{margproof}{1}}

\newcounter{margdetails}
\ndef{\margdetails}{\marginpar{\bf details} \addtocounter{margdetails}{1}}

\newcounter{margproofb}
\ndef{\margproofb}{\marginpar{\bf \small Proof (B)} \addtocounter{margproof}{1}}

\newcounter{margdetailsb}
\ndef{\margdetailsb}{\marginpar{\bf \small Details (B)} \addtocounter{margdetailsb}{1}}

\ndef{\mytimes}{\!\times\!}
\ndef{\sss}[1]{\subsubsection{}\label{#1}}
\rndef{\phi}{\varphi}
\ndef{\OpenUnitDisk}{D}
\ndef{\RHS}{RHS}                            
\ndef{\LHS}{LHS} 
\ndef{\ttt}{\Leftrightarrow}
\ndef{\then}{\Rightarrow}
\ndef{\tto}{\longrightarrow}
\ndef{\nno}{\nonumber\\}
\ndef{\newn}[1]{\index{#1} \emph{#1}}       
\ndef{\la}{\langle}
\ndef{\ra}{\rangle}
\ndef{\dbar}{{\;\bar{\phantom{o}} \!\!\!\! d}}
\ndef{\stl}[1]{\stackrel{\vbox to 0pt{\vss\hbox{$\scriptstyle #1$}}}}
\ndef{\mathcomment}[1]{{\scriptstyle\text{(#1)}}\qquad}        
\ndef{\details}[1]{\smallskip\begin{center} {\bf Here:} #1\end{center}\medskip}
\ndef{\indexcom}[1]{ --- #1}
\ndef{\longsim}{\ \sim \ }              
\ndef{\tire}{-}              
\ndef{\intinfinf}{\int_{-\infty}^\infty}
\ndef{\refnsftrace}{\cite[V.\,2.\,1]{TakI}} 
\ndef{\refaffloper}{\cite[IV.\,5, Exercise 3]{TakI}} 
\ndef{\refsemifinvNa}{\cite[V.\,1.\,21]{TakI}} 
\ndef{\reftaumeasurable}{\cite[Definition 1.2]{FK86PJM}} 
\ndef{\reftautraceclassaffl}{\cite[V.2, p.\,320]{TakI}} 
\ndef{\refinvoperideal}{\cite[Appendix A.2]{CP2}} 
\ndef{\reftautracenorm}{\cite[V.2, p.\,320]{TakI}} 
\ndef{\reftaucompact}{\cite{}} 
\ndef{\reftauFredholm}{\cite[Appendix B]{PR94JFA}} 

     \ndef{\npartial}{\slash\!\!\!\partial}
     \ndef{\Heis}{\operatorname{Heis}}
     \ndef{\Solv}{\operatorname{Solv}}
     \ndef{\Spin}{\operatorname{Spin}}
     \ndef{\SO}{\operatorname{SO}}
     \ndef{\Index}{\operatorname{index}}

             \ndef{\coker}{{\mbox coker}}
             \ndef{\p}{\partial}
             \ndef{\dd}{|\clD|}
             \ndef{\n}{\parallel}

     \ndef{\gf}[2]{\genfrac{}{}{0pt}{}{#1}{#2}}
     \ndef{\ta}{\widetilde{\alpha}}
     \ndef{\tb}{\widetilde{\beta}}
     \ndef{\txi}{\widetilde{\xi}}
     \ndef{\tk}{\widetilde{K}}
     \ndef{\CGh}{\widetilde{\CG}}
     \ndef{\boe}{{\bf e}}\ndef{\bt}{{\bf t}}
     \ndef{\vth}{\vartheta}
     \ndef{\db}{\overline{\partial}}
     \ndef{\hV}{\hat{V}}
     \ndef{\cag}{{\clA^\Gamma}}
     \ndef{\sind}{\sigma{\rm -ind}}

\let\LatexCite=\cite  

\let\ifnumref\iftrue 

\ndef{\ifuncited}[4]{\expandafter\ifx\csname used#4\endcsname\relax}

\ndef{\ifcited}[4]{\expandafter\ifx\csname used#4\endcsname\relax\else}

  \ndef{\papertitle}[1]{ \emph{#1}, }
  \ndef{\paperauthor}[2]{#2}  
  \ndef{\pbbi}[9]{%
      \ifcited{#1}{#2}{#3}{#5}%
        \ifnumref%
          \bibitem{#5}\paperauthor{#1}{#6},\papertitle{#7}#8.%
        \else%
          \advance #9 by 1%
          \ifnum#9<1%
            \bibitem[#4]{#5}\paperauthor{#1}{#6}, \papertitle{#7}#8.%
          \else%
            \bibitem[#4$\!_{\the#9}\!$]{#5}\paperauthor{#1}{#6},\papertitle{#7}#8.%
          \fi%
        \fi%
      \fi%
  }
  \ndef{\mbbi}[8]{%
     \ifcited{#1}{#2}{#3}{#5}%
        \ifnumref%
          \bibitem{#5}\paperauthor{#1}{#6},\papertitle{#7}#8.%
        \else%
          \bibitem[#4]{#5}\paperauthor{#1}{#6},\papertitle{#7}#8.%
        \fi%
     \fi%
  }

\ndef{\AddCite}[1]{%
   \ifuncited{0}{0}{0}{#1}%
     \expandafter\gdef\csname used#1\endcsname {}%
   \fi%
}

\def\ProcessCite#1,{%
     \ifx\relax#1%
         \let\next=\relax%
     \else%
         \AddCite{#1}%
         \let\next=\ProcessCite%
     \fi%
     \next%
}

\ndef{\AddCites}[1]{\ProcessCite#1,\relax,}

\ndef{\CiteWithoutExtension}[1]{%
   \AddCites{#1}%
   \LatexCite{#1}%
}

\def\CiteWithExtension[#1]#2{%
   \AddCites{#2}%
   \LatexCite[#1]{#2}%
}

\ndef{\CleverCite}{%
    \ifx\NChar[ %
       \let\MyCite=\CiteWithExtension %
    \else %
       \let\MyCite=\CiteWithoutExtension %
    \fi %
    \MyCite%
}

\renewcommand{\cite}{\futurelet\NChar\CleverCite}

      \ndef{\volume}[1]{{\bf #1}}
      \ndef{\VolYearPP}[3]{\ifnum#2=0 (to appear)\else\volume{#1} (#2), #3\fi}
      \ndef{\VolNoYearPP}[4]{\ifnum#3=0 (to appear)\else\volume{#1} #2 (#3), #4\fi}
      \ndef{\libcode}[1]{}

\ndef{\jnActaMath}[3]{Acta Math. \VolYearPP{#1}{#2}{#3}}                       
\ndef{\jnAdvMath}[3]{Adv. in~Math. \VolYearPP{#1}{#2}{#3}}                     
\ndef{\jnAlgAnal}[3]{Algebra i~Analiz \VolYearPP{#1}{#2}{#3}}
\ndef{\jnAmerMathMonth}[3]{Amer. Math. Monthly \VolYearPP{#1}{#2}{#3}}         
\ndef{\jnAnnMath}[4]{Ann. of~Math. \VolNoYearPP{#1}{#2}{#3}{#4}}               
\ndef{\jnAnalMath}[3]{J. Anal. Math. \VolYearPP{#1}{#2}{#3}}                   
\ndef{\jnBullLondMathSoc}[3]{Bull. London Math. Soc. \VolYearPP{#1}{#2}{#3}}   
\ndef{\jnBullAMS}[3]{Bull. Amer. Math. Soc. \VolYearPP{#1}{#2}{#3}}   
\ndef{\jnCanMathBull}[3]{Canad. Math. Bull. \VolYearPP{#1}{#2}{#3}}            
\ndef{\jnCanMath}[3]{Canad. J.~Math. \VolYearPP{#1}{#2}{#3}}             
\ndef{\jnCommMathPhys}[3]{Comm. Math. Phys \VolYearPP{#1}{#2}{#3}}             
\ndef{\jnCommPDE}[3]{Comm. Partial Differential Equations \VolYearPP{#1}{#2}{#3}}             
\ndef{\jnComptRendue}[3]{C.\,R.~Acad. Sci. Paris S\'er. A-B \VolYearPP{#1}{#2}{#3}}      
\ndef{\jnDiffGeom}[3]{J.~Diff. Geom. \VolYearPP{#1}{#2}{#3}}                   
\ndef{\jnErgodicTheory}[3]{Ergodic Theory and Dynamical Systems \VolYearPP{#1}{#2}{#3}} 
\ndef{\jnFuncAnal}[3]{J.~Functional Analysis \VolYearPP{#1}{#2}{#3}}           
\ndef{\jnFunkAnalPril}[4]{Функциональный анализ и его приложения \VolNoYearPP{#1}{#2}{#3}{#4}}  
\ndef{\jnGAFA}[3]{GAFA \VolYearPP{#1}{#2}{#3}}                                 
\ndef{\jnIHES}[3]{IHES Publ. Math. (Paris) \VolYearPP{#1}{#2}{#3}}             
\ndef{\jnIEOT}[3]{Integral Equations Operator Theory   \VolYearPP{#1}{#2}{#3}} 
\ndef{\jnIsrMath}[3]{Israel J.~Math. \VolYearPP{#1}{#2}{#3}}                   
\ndef{\jnKTheory}[3]{K-Theory \VolYearPP{#1}{#2}{#3}}                          
\ndef{\jnLetMathPhys}[3]{Lett. Math. Phys. \VolYearPP{#1}{#2}{#3}}             
\ndef{\jnMathAnn}[3]{Math. Ann. \VolYearPP{#1}{#2}{#3}}                        
\ndef{\jnMathAnalAppl}[3]{J.~Math. Anal. and Appl. \VolYearPP{#1}{#2}{#3}}     
\ndef{\jnMathNachr}[3]{Math. Nachr. \VolYearPP{#1}{#2}{#3}}
\ndef{\jnMathPhys}[3]{J. Math. Phys. \VolYearPP{#1}{#2}{#3}}
\ndef{\jnMathSocJap}[3]{J. Math. Soc. Japan \VolYearPP{#1}{#2}{#3}}
\ndef{\jnOperTheory}[3]{J.~Operator Theory \VolYearPP{#1}{#2}{#3}}             
\ndef{\jnPacJMath}[3]{Pacific J.~Math. \VolYearPP{#1}{#2}{#3}}                  
\ndef{\jnPositivity}[3]{Positivity \VolYearPP{#1}{#2}{#3}}
\ndef{\jnProcAmerMS}[3]{Proc. Amer. Math. Soc. \VolYearPP{#1}{#2}{#3}}         
\ndef{\jnProcCambPhilSoc}[3]{Math. Proc. Camb. Phil. Soc. \VolYearPP{#1}{#2}{#3}}
\ndef{\jnReineAngew}[3]{J.~Reine Angew. Math. \VolYearPP{#1}{#2}{#3}}          
\ndef{\jnTokyoMath}[3]{Tokyo J.~Math. \VolYearPP{#1}{#2}{#3}}
\ndef{\jnTopology}[3]{Topology \VolYearPP{#1}{#2}{#3}}
\ndef{\jnTransAmerMathSoc}[3]{Trans. Amer. Math. Soc. \VolYearPP{#1}{#2}{#3}}
\ndef{\jnIzvANSSSR}[3]{Izv. Akad. Nauk SSSR, Ser. Mat. \VolYearPP{#1}{#2}{#3}}
\ndef{\jnIzvVyshUchZav}[3]{Izv. Vyssh. Uch. Zav., Mat. \VolYearPP{#1}{#2}{#3} (Russian)}
\ndef{\jnIzdatLenUniv}[2]{Izdat. Leningrad. Univ., Leningrad, (#1), #2 (Russian)}
\ndef{\jnFieldsInsComm}[3]{Fields Inst. Comm. \VolYearPP{#1}{#2}{#3}}
\ndef{\jnDoklANSSSR}[3]{Dokl. Akad. Nauk SSSR \VolYearPP{#1}{#2}{#3}}
\ndef{\jnMatZametki}[3]{Matem. zametki \VolYearPP{#1}{#2}{#3}}
\ndef{\jnRussMathSurvey}[3]{Russian Math. Surveys \VolYearPP{#1}{#2}{#3}}
\ndef{\jnSibMathJ}[3]{Sib. Math.~J. \VolYearPP{#1}{#2}{#3}}
\ndef{\jnSovMath}[3]{J.~Soviet math. \VolYearPP{#1}{#2}{#3}}
\ndef{\jnTransMoscMathSoc}[3]{Trans. Moscow Math. Soc. \VolYearPP{#1}{#2}{#3}}
\ndef{\jnUMN}[3]{Uspekhi Mat. Nauk \VolYearPP{#1}{#2}{#3}}

\ndef{\bkTransMathMon}[2]{Trans. Math. Monographs, AMS, \volume{#1}, #2}

\ndef{\pbBirkhauser}[1]{Birkh\"auser, Boston, #1}
\ndef{\pbFactorial}[1]{Moscow, Factorial, #1}
\ndef{\pbGauthier}[1]{Gauthier-Villars, Paris, #1}
\ndef{\pbNauka}[1]{Moscow, Nauka, #1 (Russian)}
\ndef{\pbNaukaR}[1]{Москва, Наука, #1}
\ndef{\pbPrinceton}[1]{Princeton University Press, Princeton, New Jersey, #1}
\ndef{\pbPublPerish}[1]{Publish or Perish Inc., Berkeley, #1}
\ndef{\pbSpringer}[1]{Springer-Verlag, #1}

\ndef{\myauthor}[1]{\mbox{#1}}

\ndef{\Agmon}{\myauthor{Sh.\,Agmon}}
\ndef{\Ahiezer}{\myauthor{N.\,I.\,Ahiezer}}
\ndef{\Arazy}{\myauthor{J.\,Arazy}}
\ndef{\Astashkin}{\myauthor{S.\,V.\,Astashkin}}
\ndef{\Atiyah}{\myauthor{M.\,Atiyah}}
\ndef{\Avron}{\myauthor{J.\,Avron}}
\ndef{\Azamov}{\myauthor{N.\,A.\,Azamov}}
\ndef{\Banach}{\myauthor{S.\,Banach}}
\ndef{\Benameur}{\myauthor{M-T.\,Benameur}}
\ndef{\Bennett}{\myauthor{C.\,Bennett}}
\ndef{\Berezin}{\myauthor{F.\,A.\,Berezin}}
\ndef{\Berline}{\myauthor{N.\,Berline}}
\ndef{\Birman}{\myauthor{M.\,Sh.\,Birman}}
\ndef{\Blackadar}{\myauthor{B.\,Blackadar}}
\ndef{\Bogolyubov}{\myauthor{N.\,N.\,Bogolyubov}}
\ndef{\Bonsall}{\myauthor{F.\,F.\,Bonsall}}
\ndef{\BoosBavnbek}{\myauthor{B.\,Boo$\beta$-Bavnbek}}
\ndef{\Bott}{\myauthor{R.\,Bott}}
\ndef{\Bratteli}{\myauthor{O.\,Bratteli}}
\ndef{\Bredon}{\myauthor{G.\,E.\,Bredon}}
\ndef{\Breuer}{\myauthor{M.\,Breuer}}
\ndef{\Brown}{\myauthor{L.\,G.\,Brown}}
\ndef{\Bruneau}{\myauthor{V.\,Bruneau}}
\ndef{\Buslaev}{\myauthor{V.\,S.\,Buslaev}}
\ndef{\Carey}{\myauthor{A.\,L.\,Carey}}
\ndef{\CareyRW}{\myauthor{R.\,W.\,Carey}} 
\ndef{\Cartan}{\myauthor{H.\,Cartan}}
\ndef{\Chilin}{\myauthor{V.\,I.\,Chilin}}
\ndef{\Coburn}{\myauthor{L.\,A.\,Coburn}}
\ndef{\Connes}{\myauthor{A.\,Connes}}
\ndef{\Cornfeld}{\myauthor{I.\,P.\,Cornfeld}}
\ndef{\Daletskii}{\myauthor{Yu.\,L.\,Daletski\u\i}}   
\ndef{\Dixmier}{\myauthor{J.\,Dixmier}}
\ndef{\DoddsPG}{\myauthor{P.\,G.\,Dodds}}
\ndef{\DoddsTK}{\myauthor{T.\,K.\,Dodds}}
\ndef{\Douglas}{\myauthor{R.\,G.\,Douglas}}
\ndef{\Dubrovin}{\myauthor{B.\,A.\,Dubrovin}}
\ndef{\Dugundji}{\myauthor{J.\,Dugundji}}
\ndef{\Duncan}{\myauthor{J.\,Duncan}}
\ndef{\Dunford}{\myauthor{N.\,Dunford}}
\ndef{\Dykema}{\myauthor{K.\,J.\,Dykema}}
\ndef{\Edwards}{\myauthor{R.\,E.\,Edwards}}
\ndef{\Eilenberg}{\myauthor{S.\,Eilenberg}}
\ndef{\Entina}{\myauthor{S.\,B.\,\`Entina}}
\ndef{\Fack}{\myauthor{T.\,Fack}} 
\ndef{\Faddeev}{\myauthor{L.\,D.\,Faddeev}}
\ndef{\Farber}{\myauthor{M.\,Farber}}
\ndef{\Farforovskaya}{\myauthor{Yu.\,B.\,Farforovskaya}}
\ndef{\Federer}{\myauthor{H.\,Federer}}
\ndef{\Fedosov}{\myauthor{B.\,V.\,Fedosov}}
\ndef{\Figiel}{\myauthor{T.\,Figiel}} 
\ndef{\Figueroa}{\myauthor{H.\,Figueroa}}
\ndef{\Fillmore}{\myauthor{P.\,A.\,Fillmore}}
\ndef{\Fomenko}{\myauthor{A.\,T.\,Fomenko}} 
\ndef{\Fomin}{\myauthor{S.\,V.\,Fomin}}
\ndef{\Frohlich}{\myauthor{J.\,Fr\"ohlich}}
\ndef{\Fuglede}{\myauthor{B.\,Fuglede}}
\ndef{\Furutani}{\myauthor{K.\,Furutani}}
\ndef{\Gelfand}{\myauthor{I.\,M.\,Gelfand}}
\ndef{\Gesztesy}{\myauthor{F.\,Gesztesy}}     
\ndef{\Getzler}{\myauthor{E.\,Getzler}} 
\ndef{\Gilkey}{\myauthor{P.\,B.\,Gilkey}}
\ndef{\Gitler}{\myauthor{S.\,Gitler}}
\ndef{\Glazman}{\myauthor{I.\,M.\,Glazman}}
\ndef{\Glimm}{\myauthor{J.\,Glimm}}
\ndef{\Gohberg}{\myauthor{I.\,C.\,Gohberg}}
\ndef{\Golze}{\myauthor{F.\,Golze}}
\ndef{\GraciaBondia}{\myauthor{J.\,M.\,Gracia-Bond\'{i}a}}
\ndef{\Greenleaf}{\myauthor{F.\,P.\,Greenleaf}}
\ndef{\Gromov}{\myauthor{M.\,Gromov}}
\ndef{\Gunning}{\myauthor{R.\,C.\,Gunning}}
\ndef{\Haagerup}{\myauthor{U.\,Haagerup}}
\ndef{\Haag}{\myauthor{R.\,Haag}}
\ndef{\Halmos}{\myauthor{Halmos}}
\ndef{\Hardy}{\myauthor{G.\,H.\,Hardy}}
\ndef{\Higson}{\myauthor{N.\,Higson}}  
\ndef{\Hoermander}{\myauthor{L.\,Hoermander}} 
\ndef{\Hoffman}{\myauthor{K.\,Hoffman}} 
\ndef{\Ito}{\myauthor{K.\,Ito}}
\ndef{\Jaffe}{\myauthor{A.\,Jaffe}}
\ndef{\James}{\myauthor{I.\,M.\,James}}
\ndef{\Javrjan}{\myauthor{V.\,A.\,Javrjan}}
\ndef{\Kadison}{\myauthor{R.\,V.\,Kadison}}
\ndef{\Kalton}{\myauthor{N.\,J.\,Kalton}} 
\ndef{\Kato}{\myauthor{T.\,Kato}} 
\ndef{\Kobayashi}{\myauthor{S.\,Kobayashi}}
\ndef{\Koplienko}{\myauthor{L.\,S.\,Koplienko}}
\ndef{\Korotyaev}{\myauthor{E.\,Korotyaev}}
\ndef{\Kosaki}{\myauthor{H.\,Kosaki}}
\ndef{\Krein}{\myauthor{Kre\u\i n}}
\ndef{\KreinMG}{\myauthor{M.\,G.\,Kre\u\i n}}
\ndef{\KreinSG}{\myauthor{S.\,G.\,Kre\u\i n}}
\ndef{\Kuroda}{\myauthor{S.\,T.\,Kuroda}}
\ndef{\Leichtnam}{\myauthor{E.\,Leichtnam}}
\ndef{\Lesch}{\myauthor{M.\,Lesch}}
\ndef{\Lesniewski}{\myauthor{A.\,Lesniewski}}
\ndef{\Levitan}{\myauthor{B.\,M.\,Levitan}}
\ndef{\Lidskii}{\myauthor{V.\,B.\,Lidskii}}
\ndef{\Lifshits}{\myauthor{I.\,M.\,Lifshits}}
\ndef{\Lindenstrauss}{\myauthor{J.\,Lindenstrauss}}
\ndef{\Loday}{\myauthor{J.-L.\,Loday}}
\ndef{\Lord}{\myauthor{S.\,Lord}}      
\ndef{\Lorentz}{\myauthor{G.\,Lorentz}}
\ndef{\Magnus}{\myauthor{W.\,Magnus}}
\ndef{\Makarov}{\myauthor{K.\,A.\,Makarov}}
\ndef{\Mathai}{\myauthor{V.\,Mathai}}         
\ndef{\McKean}{\myauthor{H.\,P.\,McKean}}
\ndef{\Mishchenko}{\myauthor{A.\,S.\,Mishchenko}}
\ndef{\Moore}{\myauthor{C.\,C.\,Moore}}
\ndef{\Moscovici}{\myauthor{H.\,Moscovici}}  
\ndef{\Motovilov}{\myauthor{A.\,K.\,Motovilov}}
\ndef{\Moyer}{\myauthor{R.\,D.\,Moyer}}
\ndef{\Naboko}{\myauthor{S.\,N.\,Naboko}}
\ndef{\Narasimhan}{\myauthor{R.\,Narasimhan}}
\ndef{\Nomizu}{\myauthor{K.\,Nomizu}}
\ndef{\Novikov}{\myauthor{S.\,P.\,Novikov}}
\ndef{\Osterwalder}{\myauthor{K.\,Osterwalder}}
\ndef{\Patodi}{\myauthor{V.\,Patodi}}
\ndef{\Pagter}{\myauthor{B.\,de~Pagter}}  
\ndef{\Pavlov}{\myauthor{B.\,S.\,Pavlov}}
\ndef{\Pedersen}{\myauthor{G.\,K.\,Pedersen}}
\ndef{\Peller}{\myauthor{V.\,V.\,Peller}}
\ndef{\Perera}{\myauthor{V.\,S.\,Perera}}
\ndef{\Petunin}{\myauthor{Ju.\,I.\,Petunin}}
\ndef{\Phillips}{\myauthor{J.\,Phillips}}  
\ndef{\Piazza}{\myauthor{P.\,Piazza}}   
\ndef{\Pincus}{\myauthor{J.\,D.\,Pincus}}   
\ndef{\Poincare}{Poincar\'e}
\ndef{\Postnikov}{\myauthor{M.\,M.\,Postnikov}} 
\ndef{\Prinzis}{\myauthor{R.\,Prinzis}}
\ndef{\Privalov}{\myauthor{I.\,I.\,Privalov}}
\ndef{\Pushnitski}{\myauthor{A.\,B.\,Pushnitski}} 
\ndef{\Raeburn}{\myauthor{I.\,Raeburn}}
\ndef{\Raikov}{\myauthor{G.\,Raikov}}
\ndef{\Reed}{\myauthor{M.\,Reed}}
\ndef{\Rennie}{\myauthor{A.\,Rennie}}
\ndef{\Rickart}{\myauthor{C.\,E.\,Rickart}}
\ndef{\Riesz}{\myauthor{F.\,Riesz}}
\ndef{\Ringrose}{\myauthor{J.\,Ringrose}}
\ndef{\Robinson}{\myauthor{D.\,Robinson}}
\ndef{\Rossi}{\myauthor{H.\,Rossi}}
\ndef{\Rudin}{\myauthor{W.\,Rudin}}
\ndef{\Ruelle}{\myauthor{D.\,Ruelle}}
\ndef{\Ruzhansky}{\myauthor{M.\,Ruzhansky}}
\ndef{\Sakai}{\myauthor{Sh.\,Sakai}}
\ndef{\Sargsjan}{\myauthor{I.\,S.\,Sargsjan}}
\ndef{\Sato}{\myauthor{H.\,Sato}}
\ndef{\Schaeffer}{\myauthor{D.\,G.\,Schaeffer}}
\ndef{\Schluchtermann}{\myauthor{G.\,Schluchtermann}}
\ndef{\Schochet}{\myauthor{C.\,Schochet}}
\ndef{\Schrodinger}{\myauthor{E.\,Schr\"odinger}}
\ndef{\Schroodinger}{\myauthor{Schr\"odinger}}
\ndef{\Schrohe}{\myauthor{E.\,Schrohe}}
\ndef{\Schwartz}{\myauthor{J.\,T.\,Schwartz}}
\ndef{\Sedaev}{\myauthor{A.\,A.\,Sedaev}}
\ndef{\Seiler}{\myauthor{R.\,Seiler}}
\ndef{\Semenov}{\myauthor{E.\,M.\,Semenov}}
\ndef{\Shabat}{\myauthor{B.\,V.\,Shabat}}
\ndef{\Shafarevich}{\myauthor{I.\,R.\,Shafarevich}}
\ndef{\Sharpley}{\myauthor{R.\,Sharpley}}
\ndef{\Shilov}{\myauthor{G.\,E.\,Shilov}}
\ndef{\Shirkov}{\myauthor{D.\,V.\,Shirkov}}
\ndef{\Shubin}{\myauthor{M.\,A.\,Shubin}}
\ndef{\Silverman}{\myauthor{H.\,Silverman}}
\ndef{\Simon}{\myauthor{B.\,Simon}}
\ndef{\Sinai}{\myauthor{Ya.\,G.\,Sinai}}
\ndef{\Singer}{\myauthor{I.\,M.\,Singer}}
\ndef{\Solomyak}{\myauthor{M.\,Z.\,Solomyak}}
\ndef{\Soloviev}{\myauthor{Yu.\,P.\,Soloviev}}
\ndef{\Spivak}{\myauthor{M.\,Spivak}}
\ndef{\Stenkin}{\myauthor{V.\,V.\,Sten'kin}}
\ndef{\Stratila}{\myauthor{S.\,Stratila}}
\ndef{\Sucheston}{\myauthor{L.\,Sucheston}}
\ndef{\Sukochev}{\myauthor{F.\,A.\,Sukochev}}
\ndef{\Switzer}{\myauthor{R.\,M.\,Switzer}}
\ndef{\SzNagy}{\myauthor{B.\,Sz.-Nagy}}
\ndef{\Takesaki}{\myauthor{M.\,Takesaki}}
\ndef{\Taylor}{\myauthor{M.\,E.\,Taylor}}
\ndef{\Treves}{\myauthor{F.\,Treves}}
\ndef{\Troitsky}{\myauthor{E.\,V.\,Troitsky}}
\ndef{\Tzafriri}{\myauthor{L.\,Tzafriri}}
\ndef{\Varilly}{\myauthor{J.\,C.\,V\'{a}rilly}}
\ndef{\Vergne}{\myauthor{M.\,Vergne}}
\ndef{\Vladimirov}{\myauthor{V.\,S.\,Vladimirov}}
\ndef{\Voiculescu}{\myauthor{D.\,Voiculescu}}
\ndef{\Weiss}{\myauthor{G.\,Weiss}}
\ndef{\Wells}{\myauthor{R.\,O.\,Wells}}
\ndef{\Williams}{\myauthor{J.\,P.\,Williams}}
\ndef{\Winkler}{\myauthor{S.\,Winkler}}
\ndef{\Witten}{\myauthor{E.\,Witten}}
\ndef{\Wodzicki}{\myauthor{M.\,Wodzicki}}
\ndef{\Wojciechowski}{\myauthor{K.\,P.\,Wojciechowski}}
\ndef{\Yafaev}{\myauthor{D.\,R.\,Yafaev}}
\ndef{\Yosida}{\myauthor{K.\,Yosida}}
\ndef{\Zsido}{\myauthor{L.\,Zsido}}

\newcommand{\Cc}{C_\csupp^\infty(\mbR)}
\newcommand{\Texp}{\mathrm{T}\!\exp}
\newcommand{\Phia}{\Phi^{(a)}}
\newcommand{\Phis}{\Phi^{(s)}}
\newcommand{\xia}{\xi^{(a)}}
\newcommand{\xis}{\xi^{(s)}}
\newcommand{\mua}{\mu^{(a)}}
\newcommand{\mus}{\mu^{(s)}}
\newcommand{\LambdaA}{{\Lambda_\clA}}

\rndef{\LpH}[1]{\frS_{#1}\hilbasargument}       

\begin{document}
\title[Infinitesimal spectral flow]{Infinitesimal spectral flow \\ and scattering matrix}
\author{\Azamov}
\address{School of Informatics and Engineering
   \\ Flinders University of South Australia
   \\ Bedford Park, 5042, SA Australia.}
\email{azam0001@infoeng.flinders.edu.au}
\keywords{Spectral shift function, scattering matrix, infinitesimal spectral flow,
    trace compatible operators, Birman-\Krein\ formula
}
\subjclass[2000]{ 
    Primary 47A55; 
    Secondary 47A11 
}
\begin{abstract} In this note we introduce the absolutely continuous
and the singular parts of the spectral shift function as integrals
of the absolutely continuous and, respectively, of the singular parts of the infinitesimal spectral flow.
Under certain assumption, we show that this definition is independent of the piecewise linear path of integration.
The proof is based on a representation of the scattering operator of a pair of trace
compatible operators as chronological exponential of the infinitesimal
scattering matrix, and on the fact that the trace of the infinitesimal scattering matrix is equal
to the absolutely continuous part of the infinitesimal spectral flow. As a corollary,
a variant of the Birman-\Krein\ formula is derived. An interpretation of Pushnitski's
$\mu$-invariant is given.
\end{abstract}
\maketitle
\begin{center} {\small \sf\today} \end{center}

\section*{Introduction}
Let $H_0$ be a self-adjoint operator, and let $V$ be a trace class operator on a Hilbert space $\hilb.$
Then \KreinMG 's famous result \cite{Kr53MS} says that there is a unique $L^1$\tire function
$\xi_{H_0+V,H_0}(\lambda),$ known as the Lifshits-\Krein\ spectral shift function, such that for any $\Cc$\tire function $f$
\begin{gather} \label{F: Krein's formula}
  \Tr(f(H_0+V) - f(H_0)) = \int_{-\infty}^\infty f'(\lambda) \xi_{H_0+V,H_0}(\lambda)\,d\lambda.
\end{gather}
The notion of the spectral shift function was discovered by the physicist \Lifshits\ \cite{Li52UMN}.
An excellent survey on the theory of the spectral shift function can be found in \cite{BP98IEOT}.

In 1975, Birman and Solomyak \cite{BS72SM} proved the following remarkable formula for the spectral shift function
\begin{gather} \label{F: BS formula}
  \xi(\lambda) = \frac d{d\lambda} \int_0^1 \Tr(V E_{(-\infty,\lambda]}^{H_r})\,dr,
\end{gather}
where $H_r = H_0 + rV,$ $r \in \mbR,$ and $E_{(-\infty,\lambda]}^{H_r}$ is the spectral projection
(see also \cite{Si98PAMS}).
On the basis of the Birman-Solomyak formula,
the notion of infinitesimal spectral flow was introduced in \cite{AS2}.
The infinitesimal spectral flow of a self-adjoint operator $H$ under
a trace compatible perturbation $V$ (see the definition in the text)
is defined by the formula
$$
  \Phi_H(V)(\phi) = \Tr(V\phi(H)), \quad \phi \in \Cc.
$$
It was shown in \cite{AS2} that for any two operators $H_0,H_1 \in \clA$
from a trace compatible affine space $\clA$ (see the definition in the text) one can
define the spectral shift function $\xi_{H_1,H_0}$ of the pair $H_0,H_1$ as
the integral of the infinitesimal spectral flow
$$
  \xi_{H_1,H_0}(\phi) = \int_\Gamma \Phi_{H_r}(\dot H_r)(\phi)\,dr, \quad \phi \in \Cc,
$$
where $\Gamma = \set{H_r}$ is any piecewise smooth path in $\clA,$ connecting $H_0$ and $H_1.$
It was shown in \cite{AS2} that the integral does not depend on the choice of the path~$\Gamma,$
and the spectral shift function, as defined above,
is an absolutely continuous measure.
Examples of trace compatible affine spaces include
the classical case of $H_0 + \saLpH{1}$ with an arbitrary self-adjoint
operator $H_0$ on $\hilb,$ $D_0 + \saBH$ with a self-adjoint operator $D_0$ with compact resolvent \cite{ACS}
and an affine space of Schr\"odinger operators of the form $-\Delta + V + \ell^1(L^\infty),$ where $V \in L^\infty(\mbR^d)$
\cite[Section B9]{Si82BAMS} (for definition of $\ell^1(L^\infty)$ see e.g. \cite[Chapter 4]{SimTrId}).

The well-known Birman-\Krein\ formula for the spectral shift function (\cite{BK62DAN}, see also \cite[\S\,8.4]{Ya})
asserts that for a.e. $\lambda \in \mbR$
\begin{gather} \label{F: Birman-Krein formula}
  \det S(\lambda) = e^{-2\pi i \xi(\lambda)},
\end{gather}
where $S(\lambda)$ is the scattering matrix of the pair $H_0, H_0 + V,$ $V \in \LpH{1}.$
This formula was discovered for the first time by \Buslaev\ and \Faddeev\
in the case of Sturm-Liouville operators on a half-line \cite{BF60DAN}.

In this note we introduce the absolutely continuous, $\Phia,$ and
the singular, $\Phis,$ parts
of the infinitesimal spectral flow $\Phi$ by formulas (\ref{F: def of Phia}) and (\ref{F: def of Phis}),
and a decomposition of the spectral shift function
$$
  \xi_{H_1,H_0} = \xia_{H_1,H_0} + \xis_{H_1,H_0},
$$
where $\xia = \int_\Gamma \Phia$ and $\xis = \int_\Gamma \Phis,$
\ $\Gamma$ being a piecewise linear path connecting $H_0$ and $H_1.$
Similar notions in the context of Herglotz functions were considered also in \cite{GM03AA}, where one
can also find historical comments on the subject.

Under certain assumption, which includes a class of Schr\"odinger operators,
it is proved that the definition of $\xia$ and $\xis$ does not depend
on the choice of the piecewise linear path $\Gamma.$
The proof is based on the following formula for the scattering matrix
$$
  \bfS(H_1,H_0) = \Texp\brs{-2\pi i \int_0^1 W_+(H_0,H_r)\Pi_{H_r}(\dot H_r)W_+(H_r,H_0)\,dr},
$$
where $\set{H_r}_{r \in [0,1]}$ is a piecewise linear path in $\clA,$ connecting $H_0$ and $H_1.$
Though this formula is an almost straightforward consequence
of the stationary formula for the scattering matrix, it seems
to be new (to the best of the author's knowledge).

As is well-known, the Birman-\Krein\ formula (\ref{F: Birman-Krein formula})
determines the spectral shift function up to an integer-valued function.
We show that $\xi$ function in (\ref{F: Birman-Krein formula}) can be replaced by $\xia.$
This suggests that the above mentioned integer-valued function may be $\xis.$
At the end of the note an interpretation of Pushnitski's $\mu$\tire invariant is given.

I would like to thank the referee for her/his helpful comments.

\section{Notation and preliminaries}

\subsection{Trace compatibility}
We recall the notion of a trace compatible affine space of operators, which was introduced in~\cite{AS2}.
Let $\clA = H_0 + \clA_0$ be an affine space of self-adjoint operators
on a Hilbert space $\hilb,$ where~$H_0$ is a~self-adjoint operator on $\hilb$
and $\clA_0$ is a vector subspace of the real Banach space of all
bounded self-adjoint operators on $\hilb.$
We say that $\clA$ is \emph{trace compatible}, if
for all $\phi \in \Cc ,$ $V \in \clA_0$ and $H \in \clA$
\begin{gather}  \label{F: V phi(H) in L1}
  V \phi(H) \in \LpH{1},
\end{gather}
where $\LpH{1}$ is the ideal of trace class operators,
and if $\clA_0$ is endowed with a locally convex topology which coincides with
or is stronger than the uniform topology,
such that the map $(V_1, V_2) \in \clA_0^2 \mapsto V_1 \phi(H_0+V_2)$ is $\LpH{1}$\tire continuous for all
$H_0 \in \clA$ and $\phi \in \Cc.$ In particular, $\clA$ is a locally convex affine space.
We also assume, that for any $V \in \clA_0$ there exist non-negative $V_1, V_2 \in \clA_0,$
such that $V = V_1 - V_2.$ This technical condition is satisfied in all interesting examples.
We say that a pair of operators $H_0,H_1$ is \emph{trace compatible}, if they belong to some
trace compatible affine space, we say that an operator $V$ is a trace compatible perturbation
of $H,$ if the pair $H,H+V$ is trace compatible.

The infinitesimal spectral flow is a distribution valued 1-form on a trace compatible affine space of
operators $\clA,$ defined by formula
$$
  \Phi_H(V)(\phi) = \Tr(V\phi(H)), \ \ H \in \clA, \ V \in \clA_0, \ \ \phi \in \Cc.
$$

Actually, $\Phi_H(V)$ is a measure on the spectrum of $H$ \cite{AS2}.
For example,
for the trace compatible affine space $D + C_\csupp^\infty(\mbR),$ $D = \frac 1i \frac d{dx},$
if $a \in C_\csupp^\infty(\mbR),$ then for any $v \in C_\csupp^\infty(\mbR)$
$$
  \Phi_{D+a}(v) = \frac 1{2\pi} \int_\mbR v(x)\,dx \cdot \text{Lebesgue measure}.
$$
We note, that this formula is a sort of trace analogue for $\mbR$ of Connes' formula for Dixmier trace
\cite{Co88CMP,R04KT}.

\begin{lemma}
If operators $H_0$ and $H_1$ are trace compatible, then their
essential spectra coincide.
\end{lemma}
\begin{proof}
Let $U$ be an interval, which does not intersect the essential spectrum of $H_0,$
and let $\Delta \subset U$ be a segment. The projection $E^{H_0}_\Delta$ is finite-dimensional,
and, hence, it is trace class.
If $\phi$ is a smoothed indicator of $\Delta$ whose support is a subset of $U,$
then it follows from \cite[Lemma 2.1]{AS2} that $\phi(H_1) - \phi(H_0)$ is also trace class. Hence, $\phi(H_1)$
is trace class as well, which is possible only if $\Delta$ does not intersect the essential spectrum of $H_1.$
\end{proof}
\begin{lemma} \label{L: xi is int valued outside ess sp}
 If one (and hence all) operator from $\clA$ is semibounded, then
for any $H_0, H_1 \in \clA$
the spectral shift function $\xi_{H_1,H_0}$ is integer-valued
outside of their common essential spectrum.
\end{lemma}
\begin{proof} We can assume that elements of $\clA$ are semibounded from below.
Let $[a,b] \subset \mbR \setminus \sigma_{ess}(H_0)$ be a segment,
and let $f \in C^\infty$ be a decreasing function such that $f(x) = 1$ for $x<a,$ $f(x) = 0$ for $x>b.$
It follows that $f(H_0), f(H_1) \in \clF^{0,1}$ \cite[Subsection 1.5]{ACS}.
Since $H_0, H_1$ are semibounded, $f(H_1) - f(H_0)$ is compact (even trace class)
by \cite[Lemma 2.1]{AS2}. Hence, by \cite[Theorem 3.18]{ACS},
the spectral shift function $\xi_f$ of the pair $f(H_0), f(H_1)$ is integer-valued.
The invariance principle (see e.g. \cite[(1.9)]{Pu01FA} or \cite[(13)]{AS2})
now implies that $\xi_{H_1,H_0}$ is integer valued on $[a,b].$
\end{proof}
If it is known that $E^{H_1}_\lambda - E^{H_0}_\lambda$ is compact for $\lambda \notin \sigma_{ess}(H_0),$
then the proof of this lemma shows that the spectral shift function is integer valued outside the essential spectrum.
That $E^{H_1}_\lambda - E^{H_0}_\lambda$ is compact for $\lambda \notin \sigma_{ess}(H_0)$
presumably can be shown by methods of \cite{Pu07arx}.

\subsection{Scattering theory}
We recall some notions of the mathematical scattering theory from \cite{Ya} (see also \cite{BW,BYa92AA2,RS3}).
The \emph{wave operators} of a pair of self-adjoint operators $H_0, H$
on $\hilb$ are the operators \cite[(2.1.1)]{Ya}
(when the respective limits exist)
$$
  W_{\pm}(H,H_0) := \slim_{t \to \pm \infty} e^{itH} e^{-itH_0}P_a(H_0),
$$
where $P_a(H_0)$ denotes the projection onto the absolutely continuous subspace $\hilba_0$ of $H_0,$
and the limit is the strong operator limit. Further,
$W_{\pm}^*(H,H_0) = W_{\pm}(H_0,H)$ \cite[(2.2.2)]{Ya}.
The wave operator $W_{\pm}(H,H_0)$ is an isometry
of $\hilba_0$ into $\hilba,$ where $\hilba$ is the absolutely continuous subspace of $H,$
i.e.
$$
  W_{\pm}^*(H,H_0)W_{\pm}(H,H_0) = P_a(H_0) \ \ \text{and} \ \ W_{\pm}(H,H_0)W_{\pm}^*(H,H_0) \leq P_a(H).
$$
If $W_{\pm}(H,H_0)W_{\pm}^*(H,H_0) = P_a(H)$ then the wave operator $W_{\pm}(H,H_0)$
is called \emph{complete}. The completeness of $W_{\pm}(H,H_0)$ is equivalent
to the existence of $W_{\pm}(H_0,H)$ \cite[Theorem 2.3.6]{Ya}.
If the wave operators $W_\pm(H,H_0)$ exist and are complete then
the \emph{scattering operator} of a pair of operators $H_0,H$ is defined as \cite[(2.4.1)]{Ya}
$$
  \bfS(H,H_0) = W_+^*(H,H_0)W_-(H,H_0)
$$
and it is a unitary operator on the absolutely continuous subspace of $H_0.$
If the wave operators $W_\pm(H_2,H_1)$ and $W_\pm(H_1,H_0)$ exist, then
$W_\pm(H_2,H_0)$ also exists and $W_\pm(H_2,H_0) = W_\pm(H_2,H_1)W_\pm(H_1,H_0)$
\cite[Theorem 2.1.7]{Ya}.

Let $\clA$ be a trace compatible affine space of operators on a Hilbert space $\hilb.$
If $H_0, H_1 \in \clA$ then the wave operators $W_\pm(H_1,H_0)$ exist
and are complete.
Existence of $W_\pm(H_1,H_0)$ follows from Birman's local criteria
for existence of wave operators \cite[Theorem 6.4.1]{Ya}. Completeness of $W_\pm(H_1,H_0)$
follows from existence of $W_\pm(H_0,H_1).$
Hence, for any $H_0, H_1 \in \clA,$ the scattering operator $\bfS(H_1,H_0)$
exists and is a unitary operator on $\hilba_0.$
This implies that the absolutely continuous spectrum of all $H \in \clA$
coincide.

If $V$ is a trace compatible perturbation of $H,$ then
$G := \sqrt{\abs{V}}$ is $H$\tire weakly smooth \cite[Definition 5.1.1]{Ya}.
This follows from $H$\tire weak smoothness of the Hilbert-Schmidt
operator $GE^H_{\Delta}$ and \cite[Lemmas 5.1.3, 5.1.4]{Ya}.

\subsection{Direct integrals of Hilbert spaces}
Let $\clA$ be a trace compatible affine space of operators on $\hilb.$
Let $H_0 \in \clA$ and let $V \in \clA_0.$
Let $\spectruma{H_0}$ be the spectrum of the absolutely continuous part of $H_0.$
Let
\begin{gather} \label{F: direct integral of hilbert spaces}
  \euF_0 \colon \hilba_0 \to \int^{\oplus}_{\spectruma{H_0}} \hilb_\lambda\,d\lambda,
  \quad \euF_0 \colon \eta \mapsto \tilde \eta(\cdot),
\end{gather}
be an isomorphism of the Hilbert space $\hilba_0$ onto a direct integral of Hilbert spaces $\hilb_\lambda,$ such that
$(\euF_0 H_0\eta)(\lambda) = \lambda \euF_0 \eta(\lambda) = \lambda \tilde \eta(\lambda), \ \eta \in \hilba_0$ \cite[\S\,1.5]{Ya}.
The scattering operator $\bfS(H_1,H_0)$ is diagonal in this representation \cite[\S 2.4]{Ya}
and it acts as multiplication by a unitary operator $S(\lambda; H_1, H_0)$ on $\hilb_\lambda.$
The function $\lambda \mapsto S(\lambda; H_1, H_0)$ is called the \emph{scattering matrix}.

Let $\clK$ be an auxiliary Hilbert space.
Let $V = G^*JG,$ where $J = J^*$ is an invertible operator on $\clK$ and $G \colon \hilb \to \clK,$ and let
$$
  T_r(z) = G R_z(H_r)G^*, \quad z \in \mbC \setminus \mbR,
$$
where $H_0 \in \clA,$ $V \in \clA_0,$ $H_r = H_0+rV.$ One can take $\clK = \hilb,$ $G = \sqrt{\abs{V}}$ and $J = \sgn(V),$
but other decompositions are also useful \cite{Pu01FA}.

We denote by $\Lambda = \LambdaA$ the inner part of the common absolutely continuous spectrum of operators $H \in \clA.$
\begin{assump}\label{A: assumption 2}
(i) There exists $p \in [1,\infty]$ such that
for any $H \in \clA$ and $V \in \clA_0$
the function $T(z) = G R_z(H)G^*$
takes values in $\frS_p(\clK)$ for $\Im z \neq 0,$ and
for all $\lambda \in \LambdaA$ it has non-tangential limit
values $T(\lambda \pm i0)$ in $\frS_p(\clK).$
\\ (ii) The imaginary part
$B(\lambda \pm i0) := \Im T(\lambda \pm i0)$ of $T(\lambda \pm i0)$
belongs to $\frS_1(\clK)$ for all $\lambda \in \LambdaA.$
\end{assump}
This assumption holds, for example, for the space
\begin{gather} \label{F: example of clA}
  -\Delta + \set{V \in L^\infty(\mbR^d) \colon \exists \eps>0 \ \exists C > 0 \ \forall x \in \mbR^d \ \abs{V(x)} \leq C(1+\abs{x})^{-d-\eps}}
\end{gather}
with $\LambdaA = (0,\infty),$
where $\Delta$ is the Laplace operator on $\mbR^d$ with $\dom(\Delta) = \euH_2(\mbR^d)$ Sobolev space
(see e.g. \cite{BYa92AA2}, see also \cite{Agm,BeShu,Hoer2,KuJMSJ73II}).


If Assumption \ref{A: assumption 2} holds, then the formula
$$
  Z(\lambda; G)\eta = (\euF_0 G^* \eta)(\lambda), \quad  \lambda \in \LambdaA,
$$
defines operators $Z(\lambda; G) \in \clB(\clK,\hilb_\lambda)$ unambiguously \cite{Agm,BYa92AA2,Ya}.
Letting $Z_r(\lambda) = Z(\lambda; G)$ with respect to the direct integral decomposition of $H_r,$
by \cite[(5.4.4)]{Ya}
\begin{gather} \label{F: (2.7) of BYa}
  \pi Z_r^*(\lambda)Z_r(\lambda) = B_r(\lambda + i0).
\end{gather}
It follows that $Z_r(\lambda)$ is a Hilbert-Schmidt operator.
Under Assumption \ref{A: assumption 2},
if $H_0 \in \clA$ and $V \in \clA_0,$
then $V$ is an integral operator with respect to
the representation (\ref{F: direct integral of hilbert spaces}),
and has kernel
\begin{gather} \label{F: kernel of V}
  v(\lambda,\lambda') = Z_0(\lambda) J Z_0^*(\lambda'), \quad \lambda, \lambda' \in \LambdaA,
\end{gather}
which takes trace class values. This follows from \cite[Lemma 5.4.3]{Ya}
and (\ref{F: (2.7) of BYa}).

We give the proof of the following lemma for completeness (see e.g. \cite{Agm,BYa92AA2,Ya}).
\begin{lemma} \label{L: Br exists iff ...}
  If $T_0(\lambda+i0)$ exists then $T_r(\lambda+i0)$ exists if and only if the
  operator $1 + rT_0(\lambda+i0)J$ is invertible.
\end{lemma}
\begin{proof} The second resolvent identity $R(z) = R_0(z) - R(z)VR_0(z)$ implies that
$$
  T_r(z)(1+rJT_0(z)) = T_0(z).
$$
Taking limits, one gets
  \begin{gather} \label{F: T(1+T0)=T0}
    T_r(\lambda+i0)(1+rJT_0(\lambda+i0)) = T_0(\lambda+i0).
  \end{gather}
Since $T_0$ is compact, $1+rJT_0(\lambda+i0)$ is not invertible if and only if
there exists a non zero $\psi \in \hilb,$ such that $(1+rJT_0(\lambda+i0))\psi = 0.$
This and (\ref{F: T(1+T0)=T0}) imply that $T_0(\lambda+i0)\psi = 0.$ Hence $\psi = 0.$
\end{proof}

%

\section{Results}
\begin{defn}
The \emph{infinitesimal scattering matrix} $\Pi_{H_0}(V)$
of the operator $H_0 \in \clA$ under a trace compatible perturbation by $V \in \clA_0 $ is
a self-adjoint operator, which in the representation (\ref{F: direct integral of hilbert spaces})
acts by multiplication by
\begin{gather} \label{F: Pi = ZJZ*}
  \Pi_{H_0}(V)(\lambda) := v(\lambda,\lambda)
        = Z_0(\lambda)JZ_0^*(\lambda), \quad \lambda \in \LambdaA.
\end{gather}
\end{defn}

The map $V \in \clA_0 \mapsto \Pi_H(V)$ is not linear, but if perturbations $V_1, V_2 \in \clA_0$ are disjoint
(i.e. $V_1V_2 = V_2V_1 = 0$) then
\begin{gather} \label{F: Pi(V) = Pi(V1) + Pi(V2)}
  \Pi_H(V_1+V_2) = \Pi_H(V_1) + \Pi_H(V_2).
\end{gather}
If Assumption \ref{A: assumption 2} holds, then for each $\lambda \in \LambdaA$ the operator $\Pi_H(V)(\lambda)$ is trace class.

The following theorem is a classical stationary representation for the scattering matrix \cite{BE,BYa92AA2},
\cite[Theorem 5.7.1$'$]{Ya}. 
\begin{thm} If $H_0 \in \clA,$ $V \in \clA_0,$ $H_r = H_0 + rV,$ $r \in [0,1]$
and if Assumption \ref{A: assumption 2}
holds, then for all $\lambda \in \LambdaA$
and for all $r \in [0,1]$ the scattering matrix $S(\lambda; H_r,H_0)$ exists,
and the stationary representation for the scattering matrix
  \begin{gather} \label{F: stationary rep-n for SM}
    S(\lambda; H_r,H_0) = 1_\lambda - 2\pi i r Z_0(\lambda) J (1 + r T_0(\lambda+i0)J)^{-1}Z_0^*(\lambda)
  \end{gather}
holds.  Moreover, $S(\lambda) - 1_\lambda \in \frS_1(\hilb_\lambda)$
for all $\lambda \in \LambdaA.$
\end{thm}
%

\begin{lemma} \label{L: dot S = -2 pi i ...}
If $H_r = H_0 + rV,$ where $H_0 \in  \clA,$ $V \in \clA_0,$
and Assumption \ref{A: assumption 2} holds, then
$$
  \frac d{dr} S(\lambda; H_r,H_{r_0}) \big|_{r=r_0} = -2 \pi i \Pi_{H_{r_0}}(V)(\lambda)
$$
for all $\lambda \in \LambdaA,$
where the derivative is taken in $\frS_1(\hilb_\lambda)$\tire topology.
\end{lemma}
\begin{proof}
The formula (\ref{F: stationary rep-n for SM}) implies that for $\lambda \in \LambdaA$
\begin{gather} \label{F: S = 1 - 2 pi i ...}
  S(\lambda; H_r,H_{r_0}) = 1_\lambda - 2\pi i (r-r_0)Z_{r_0}(\lambda) J (1 + (r-r_0)T_{r_0}(\lambda+i0)J)^{-1}Z_{r_0}^*(\lambda),
\end{gather}
where $1+(r-r_0)T_{r_0}(\lambda+i0)J$ is invertible by Lemma \ref{L: Br exists iff ...}.
It follows that in $\frS_1(\hilb_\lambda)$
\begin{multline} \label{F: dot Sr}
  \frac d{dr} S(\lambda; H_r,H_{r_0})
     = - 2\pi i \big[Z_{r_0}(\lambda)J (1 + (r-r_0)T_{r_0}(\lambda+i0)J)^{-1}Z_{r_0}^*(\lambda)
     \\ - (r-r_0) Z_{r_0}(\lambda)J (1 + (r-r_0)T_{r_0}(\lambda+i0)J)^{-2}T_{r_0}(\lambda+i0)J Z_{r_0}^*(\lambda)\big]
   \\   = - 2\pi i \sqbrs{Z_{r_0}(\lambda)J (1 + (r-r_0)T_{r_0}(\lambda+i0)J)^{-2}Z_{r_0}^*(\lambda)}.
\end{multline}
This and (\ref{F: Pi = ZJZ*}) complete the proof.
\end{proof}

%
%
%
\begin{prop} \label{P: S = Texp}
If $H_0 \in \clA,$ $V \in \clA_0,$ if $H_r = H_0 + rV,$ $r \in [0,1]$
and if Assumption \ref{A: assumption 2} holds,
then for all $\lambda \in \LambdaA$
\begin{gather} \label{F: S = T exp...}
  S(\lambda; H_1, H_0) = \Texp\brs{-2\pi i \int_0^1 w_+(\lambda; H_0, H_r)\Pi_{H_r}(\dot H_r)(\lambda)w_+(\lambda; H_r, H_0)\,dr}.
\end{gather}
\end{prop}
\begin{proof}
It follows from Lemma \ref{L: dot S = -2 pi i ...} and \cite[Corollary 7.1.2]{Ya} that for all $\lambda \in \LambdaA$
\begin{gather*}
  \frac d{dr} S(\lambda; H_r,H_0)
    = -2\pi i w_+(\lambda; H_0, H_r) \Pi_{H_r}(V)(\lambda) w_+(\lambda; H_r,H_0)S(\lambda; H_r,H_0),
\end{gather*}
where the derivative is taken in $\frS_1(\hilb_\lambda)$\tire topology.
The derivative $\frac d{dr} S(\lambda; H_r,H_0)$ is $\frS_1(\hilb_\lambda)$\tire continuous by (\ref{F: dot Sr}).
Since $S(\lambda; H_r,H_0)$ is also $\frS_1(\hilb_\lambda)$\tire continuous, by the last formula the function
$
  r \mapsto w_+(\lambda; H_0, H_r) \Pi_{H_r}(V)(\lambda) w_+(\lambda; H_r, H_0)
$
is $\frS_1(\hilb_\lambda)$\tire continuous. Hence, integration of the last equation by Lemma \ref{L: T exp} gives (\ref{F: S = T exp...}).
\end{proof}

\begin{thm} \label{T: S = Texp for piecewise linear path}
Let Assumption \ref{A: assumption 2} holds for $\clA.$
If $H_0,H_1 \in \clA$ and if $\set{H_r}_{r \in [0,1]}$ is a piecewise linear path in $\clA$
connecting $H_0$ and $H_1,$ then for all $\lambda \in \LambdaA$ (\ref{F: S = T exp...}) holds.
\end{thm}
\begin{proof}
This follows from \cite[Corollary 7.1.2]{Ya}, Proposition \ref{P: S = Texp},
the multiplicative property of wave operators
and Lemma \ref{L: Texp su = Texp st + Texp tu}.
\end{proof}
The question of whether the piecewise linear path in this theorem
can be replaced by piecewise smooth path is open.

It is known that the infinitesimal spectral flow $\Phi_H(V)$ is a measure on $\spectrum{H}.$
Let $\Phia_H(V)$ and $\Phis_H(V)$ denote the absolutely continuous and singular parts of this measure.
In other words
\begin{gather} \label{F: def of Phia}
  \Phia_H(V)(\phi) = \Tr(V \phi(H^{(a)})), \quad \phi \in C_\csupp(\mbR),
\end{gather}
and
\begin{gather} \label{F: def of Phis}
  \Phis_H(V)(\phi) = \Tr(V \phi(H^{(s)})), \quad \phi \in C_\csupp(\mbR),
\end{gather}
where $H^{(a)}$ and $H^{(s)}$ are absolutely continuous and singular parts of $H.$
We define the "absolutely continuous" and "singular" parts of the spectral shift function $\xi$
by formulas
\begin{gather} \label{F: def of xi(a) and xi(s)}
  \xia(\phi) = \int_\Gamma \Phia(\phi), \quad \xis(\phi) = \int_\Gamma \Phis(\phi),
\end{gather}
where $\Gamma$ is any piecewise linear path in $\clA,$ connecting $H_0$ and $H_1.$ Independence of this definition
from the choice of the path $\Gamma$ will be shown in Corollary \ref{C: def of xia is good}.
The absolute continuity of $\Phia_H(V)$ follows from the following proposition.

%
%
%
\begin{prop} \label{P: Tr Pi = Phi}
If $V$ is a trace compatible perturbation of $H$ satisfying
Assumption \ref{A: assumption 2}, then for any $\phi \in \Cc$ the operator
$V\phi(H)$ is an integral operator with trace class valued kernel $v(\lambda,\lambda'),$
$\lambda,\lambda' \in \LambdaA,$ and
$$
  \Phia_H(V) (\phi) =
  \int_\LambdaA \Tr_{\hilb_\lambda}\brs{\Pi_H(V)(\lambda)} \phi(\lambda)\,d\lambda.
$$
\end{prop}
The proof of this proposition is standard (see e.g. \cite{BE,Ya}).

As a corollary we get a variant of the Birman-\Krein\ formula.
\begin{thm} \label{T: det S = exp(xi ac)}
If $H_0,H_1 \in \clA$ and Assumption \ref{A: assumption 2} holds, then
\begin{gather} \label{F: Birman Krein}
  -2\pi i \xia_{H_1,H_0}(\lambda) = \ln \det S(\lambda; H_1, H_0), \quad a.\,e.\, \ \lambda \in \LambdaA,
\end{gather}
where the branch of the logarithm is chosen in such a way, that the function
$r \in [0,1] \mapsto \ln \det S(\lambda; H_r,H_0)$ is continuous, and $\xia$ is defined
by (\ref{F: def of xi(a) and xi(s)}) with $\Gamma$ a straight line.
\end{thm}
\begin{proof}
   Let, as usual, $H_r = H_0 + rV,$ $H_0 \in \clA,$ $V \in \clA_0.$
By definition, for any $\phi \in \Cc$
\begin{gather*}
  \xia(\phi) = \int_0^1 \Phia_{H_r}(V)(\phi)\,dr
     = \int_0^1 \Tr\brs{V\phi(H^{(a)}_r)}\,dr.
\end{gather*}
Hence, by Proposition \ref{P: Tr Pi = Phi}, it follows that
$$
  \xia(\phi) = \int_0^1 \int_{\LambdaA} \Tr_{\hilb_\lambda^{(r)}}\sqbrs{v_r(\lambda,\lambda)}\phi(\lambda)\,d\lambda\,dr,
$$
where the Hilbert spaces $\hilb_\lambda^{(r)}$ are from the direct integral
decomposition (\ref{F: direct integral of hilbert spaces}) for~$H_r.$
Since the operators $w_{\pm}(\lambda; H_r,H_0) \colon \hilb_\lambda \to \hilb_\lambda^{(r)}$
are unitary \cite[Proposition 5.7.3]{Ya},
from this and Fubini's theorem it follows that
$$
  \xia(\phi) = \int_{\LambdaA} \int_0^1 \Tr_{\hilb_\lambda}\sqbrs{w_+(\lambda; H_0, H_r)v_r(\lambda,\lambda)w_+(\lambda; H_r, H_0)}
  \phi(\lambda)\,dr \,d\lambda.
$$
Proposition \ref{P: S = Texp} and Lemma \ref{L: det Texp = exp Tr} now imply 
$$
  -2\pi i \xia(\phi) = \int_{\LambdaA} \ln \det S(\lambda; H_1,H_0) \phi(\lambda) \,d\lambda,
$$
where the branch of the logarithm is chosen as in the statement of the theorem.
Since $\xi$ is absolutely continuous by \cite[Theorem 2.9]{AS2},
so is $\xia.$ Hence, for a.e. $\lambda \in \LambdaA$
$$
  -2\pi i \xia(\lambda) = \ln \det S(\lambda; H_1, H_0).
$$
\end{proof}

\begin{cor} \label{C: def of xia is good} Definition (\ref{F: def of xi(a) and xi(s)})
of $\xia$ and $\xis$ is independent of the choice of the piecewise linear path $\Gamma$ in $\clA.$
\end{cor}
\begin{proof} For $\xia$ this follows from Theorems \ref{T: S = Texp for piecewise linear path}
and \ref{T: det S = exp(xi ac)} (this also follows from the proof of Theorem \ref{T: det S = exp(xi ac)}
with reference to Theorem \ref{T: S = Texp for piecewise linear path} instead of Proposition \ref{P: S = Texp}).
For $\xis$ this follows from
$
  \xi = \xia + \xis
$ and \cite[Theorem 2.9]{AS2}.
\end{proof}

\begin{cor} For the affine space (\ref{F: example of clA})
the Birman-Krein formula (\ref{F: Birman-Krein formula}) holds.
\end{cor}
\begin{proof} Kato's theorem implies that \Schroodinger\ operators of the class (\ref{F: example of clA})
don't have singular spectrum on the absolutely continuous spectrum.
This and Theorem \ref{T: det S = exp(xi ac)} imply that for $\lambda \in \LambdaA$
(\ref{F: Birman-Krein formula}) holds.
For $\lambda \notin \LambdaA$
this follows from the equality $\sigma_{ess} = \sigma_{ac},$ which is known to be true for \Schroodinger\ operators
of the class (\ref{F: example of clA}),
and from Lemma \ref{L: xi is int valued outside ess sp}.
\end{proof}

Let $e^{i\theta_1(r)},$ $e^{i\theta_2(r)}, \ldots$ be the set of eigenvalues of $S(\lambda; H_r, H_0).$
The functions $\theta_j(r), \ j = 1,2,\ldots,$ are continuous functions of $r$ (see e.g. \cite{Ya} or \cite[IV.3.5]{Kato}).
For any $\theta \in [0,2\pi),$ one can define
Pushnitski's invariant $\mu(\theta; \lambda) = \mu(\theta; \lambda, H_1, H_0)$ as the spectral flow through
the point $e^{i\theta}$ by eigenvalues $e^{i\theta_1(r)}, e^{i\theta_2(r)}, \ldots$
of the path $\set{S(\lambda; H_r,H_0)}_{r \in [0,1]}$ \cite{Pu01FA}, i.e.
\begin{gather} \label{F: def of mu}
  \mu(\theta; \lambda) = \sum_{j=1}^\infty \brs{1 + \sqbrs{\frac{\theta_j(1) - \theta}{2\pi}}},
\end{gather}
where $[x]$ is the integer part of $x.$
The formulas (\ref{F: Birman Krein}) and (\ref{F: def of mu})
imply that for a.e. $\lambda \in \LambdaA$
\begin{gather} \label{F: xia = sum of theta's}
  \xia(\lambda) = - \frac 1{2\pi} \sum_{j=1}^\infty \theta_j(1) =  - \frac 1{2\pi } \int_0^{2\pi} \mu(\theta; \lambda)\,d\theta,
\end{gather}
which is Pushnitski's formula \cite[(1.12)]{Pu01FA}.
Theorem \ref{T: S = Texp for piecewise linear path} and \cite[\S 7.8]{Ya} imply that
the definition of $\mu(\theta; \lambda),$ given by (\ref{F: def of mu}), does not depend on the choice
of the piecewise linear path $H_r,$ connecting $H_0$ and $H_1,$ and it is well-defined
in the sense that the eigenvalues $\theta_j(r)$ do not make excessive windings around the unit circle.
We note that in \cite{Pu01FA} the scattering matrix is connected with $1$ by sending the imaginary part $y$ of
the spectral parameter $\lambda + iy$ to $+\infty.$

In \cite{Az2} it will be shown that in the case of a class of \Schroodinger\ operators which admit
embedded eigenvalues, the "singular" part $\xis$ of the spectral shift function is an
integer-valued function, which is non-zero on the absolutely continuous spectrum.
It will be shown that in this case Pushnitski's $\mu$\tire invariant admits a natural
decomposition $\mu(\theta,\lambda) = \mua(\theta,\lambda) + \mus(\theta,\lambda),$
where $\mua$ and $\mus$ are integer-valued functions, $\mua$ is non-zero only on absolutely continuous spectrum,
while $\mus$ does not depend on $\theta,$ and actually is equal to $-\xis.$
For $\mua$ the formulas (\ref{F: def of mu}) and (\ref{F: xia = sum of theta's}) hold.
In the case, considered in this note,
$\mus = \xis = 0$ on the absolutely continuous spectrum.
Proofs of these results will be based on the analysis of the behaviour of the eigenvalues
the unitary-valued functions $M(z; H_r,H_0),$ $S(z; H_0, G, rJ),$ introduced by \Pushnitski\ in \cite{Pu01FA},
and that of $S(\lambda; H_r, H_0),$ as $\Im z \to 0+$ and $\Im z \to +\infty.$
One of the elements of the proof is that the scattering matrix $S(\lambda; H_r,H_0),$ considered as a function of the coupling constant $r,$
is a meromorphic function, which admits analytical continuation to the real poles $r_0$ of $\brs{J^{-1}+rT_0(\lambda+i0)}^{-1}.$
As is well known, these poles correspond to embedded eigenvalues (see e.g. \cite[Lemma 4.7.8]{Ya} or the proof of \cite[Theorem 4.2]{Agm}).
This result can be interpreted as a jump by an integer multiple of $2\pi$
of one of the scattering phases $\theta_j(r,\lambda),$
when $r$ crosses the "resonance" point $r_0,$ i.e. the point, for which the equation $H_r \psi = \lambda \psi, \ \lambda \in \Lambda,$
has an $L^2$ solution, in accordance with Pushnitski's formula \cite[(1.12)]{Pu01FA}.
This also agrees with physical interpretation, given in \cite[XVIII.6]{Bohm}.

It seems to be likely that $\xis$ is always an integer-valued function.

\appendix

\section{Chronological exponential}
In this appendix an exposition of the chronological exponential is given. See e.g.
\cite{Gamk} and \cite[Chapter 4]{BSh}.

Let $p \in [1,\infty]$ and let $a < b.$
Let ${A(\cdot) \colon [a,b] \to \LpH{p}}$ be a piecewise continuous path of self-adjoint operators from $\LpH{p}.$
Consider the equation
\begin{gather} \label{F: dot Xt = At Xt}
  \frac {d X(t)}{dt} = \frac 1i A(t) X(t), \ \ X(a) = 1,
\end{gather}
where the derivative is taken in $\LpH{p}.$
By definition, the left chronological exponent is
\begin{gather} \label{F: def of Texp}
  \Texp\brs{\frac 1i\int_a^t A(s)\,ds} =
     1 + \sum_{k=1}^\infty \frac 1{i^k} \int_a^{t} dt_1 \int_a^{t_1} dt_2 \ldots \int_a^{t_{k-1}} dt_k A(t_1)\ldots A(t_k),
\end{gather}
where the series converges in $\LpH{p}$\tire norm.
\begin{lemma} \label{L: T exp}
The equation (\ref{F: dot Xt = At Xt}) has a unique continuous solution $X(t),$
given by formula
$$
  X(t) = \Texp\brs{\frac 1i\int_a^t A(s)\,ds}.
$$
\end{lemma}
\begin{proof} Substitution shows that (\ref{F: def of Texp}) is a continuous solution of (\ref{F: dot Xt = At Xt}).
Let $Y(t)$ be another continuous solution of (\ref{F: dot Xt = At Xt}).
Taking the integral of (\ref{F: dot Xt = At Xt}) one gets
$$
  Y(t) = 1 + \frac{1}{i} \int _a^t A(s) Y(s)\,ds.
$$
Iteration of this integral and the bound $\sup_{t \in [a,b]} \norm{A(t)}_1 \leq \const$
show that $Y(t)$ coincides with $(\ref{F: def of Texp}).$
\end{proof}
\begin{lemma} \label{L: Texp su = Texp st + Texp tu}
The following equality holds
$$
  \Texp\brs{\int_s^u A(s)\,ds} = \Texp\brs{\int_t^u A(s)\,ds} \Texp\brs{\int_s^t A(s)\,ds}.
$$
\end{lemma}
\begin{proof} Both sides of this equality are solutions of the equation
$\frac {dX(u)}{du} = \frac 1i A(u)X(u)$ with the initial condition $X(t) = \Texp\brs{\int_s^t A(s)\,ds}.$
\end{proof}
By $\det$ we denote the classical Fredholm determinant (see e.g. \cite{SimTrId}).
\begin{lemma} \label{L: det Texp = exp Tr}
If $p = 1$ then the following equality holds
$$
  \det\, \Texp\brs{\frac 1i \int_a^t A(s)\,ds} = \exp\brs{\frac 1i \int_a^t \Tr(A(s))\,ds}.
$$
\end{lemma}
\begin{proof} Let $F(t)$ and $G(t)$ be the left and the right hand sides of this equality respectively.
Then $\frac {G(t)}{dt} = \frac 1i \Tr(A(t)) G(t),$ $G(a) = 1.$
Further, by Lemma \ref{L: Texp su = Texp st + Texp tu}
$$
  \frac d{dt} F(t) = \lim_{h\to 0} \frac 1h \brs{\det \, \Texp \brs{\frac 1i \int_t^{t+h} A(s)\,ds} - 1} F(t)
  = \frac 1i \Tr(A(t)) F(t),
$$
where the last equality follows from definitions of determinant \cite[(3.5)]{SimTrId},
$\Texp$ and piecewise continuity of $A(s).$
\end{proof}


\mathsurround 0pt
\ndef{\AndSoOn}{$\dots$}


\end{document}